\documentclass{ifacconf}

\usepackage[english]{babel}
\usepackage{graphicx}      
\usepackage{natbib}        
\usepackage{subfigure}
\usepackage{color}
\usepackage{amssymb,amsmath,color}
\usepackage{srcltx,euscript}

\newcommand{\lt}{\left}
\newcommand{\rt}{\right}

\newcommand{\gs}{\geqslant}
\newcommand{\ilt}{\int\limits}

\begin{document}
\begin{frontmatter}
\title{Nonparametric gamma kernel estimates of density derivatives on positive semi-axis}
\author[Dobrovidov]{Alexander V. Dobrovidov,}
\author[Markovich]{Liubov A. Markovich}
\address[Dobrovidov]{Institute of Control Sciences, Russian Academy of Sciences, Moscow, Russia (e-mail: dobrovid@ipu.ru).}
\address[Markovich]{Institute of Control Sciences, Russian Academy of Sciences, Moscow, Russia (e-mail: kimo1@mail.ru).}
\begin{abstract}:  $\:$We consider nonparametric estimation
 of the derivative of a probability density function with the bounded support on
 $[0,\infty)$.
Estimates are looked up in the class of estimates with asymmetric
gamma kernel functions. The use of gamma kernels is due to the fact
they are nonnegative, change their shape depending on the position
on the semi-axis and possess other good properties. We found
analytical expressions for bias, variance, mean integrated squared
error (MISE) of the derivative estimate.  An optimal bandwidth, the
optimal MISE, and rate of  mean square convergence of the estimates
for density derivative have also been found.
\end{abstract}
\end{frontmatter}

\section{Introduction}
In many models of financial and actuary mathematics variables can be
only positive. That is why the proposal of adequate methods for
estimating characteristics of these models is up to date. In the
paper of Song Xi \cite{Chen:20} nonparametric gamma kernel estimate
for a reconstruction of probability density functions $f(x)$ with
support $[0,\infty)$  was proposed. As it is known, for instance,
from \cite{Jones:95} classical estimation methods with symmetric
kernels yield a large bias on the zero boundary that leads to a bad
quality of classical estimates in this case. In contrast to it,
nonparametric  estimates with asymmetric gamma kernel have a small
bias on the boundary near zero and have a variance at a point $x$ of
order $O(n^{-4/5}x^{-1/2})$, which decreases with the increase of
the argument   $x$. Such good properties of the estimates induce one
to use them as a basis for the synthesis and investigation of the
density derivative estimates.  One of the important areas of
application for density derivative estimates is the theory of
nonparametric signal estimation, published in \cite{Dobrovidov:12},
where these results are finely used, for example, in multiplicative
stochastic models. Equation for  the optimal signal estimate  are
expressed in terms of the logarithmic  density derivative of the
observed random variables which is known to contain a density
derivative. Thus, the construction and investigation of a
nonparametric kernel estimate  of the density derivative function is
the goal of this work.
\section{Main results}
Let $X_1...X_n$ be a sample of i.i.d random variables from a
distribution with an unknown probability density function  $f(x)$,
which is defined on the  support $x\in[0,\infty)$. The gamma kernel
estimate  is defined in Song Xi \cite{Chen:20} as
 \vspace{-5mm}
\begin{eqnarray}\label{1}
\hat{f}(x) = \frac{1}{n}\sum_{i=1}^{n}K_{\rho_b(x),b}(X_i),
\end{eqnarray}
where \vspace{-10mm}
\begin{eqnarray*}K_{\rho_b(x),b}(t)=\frac{t^{\rho_b(x)-1}\exp(-t/b)}{b^{\rho_b(x)}\Gamma(\rho_b(x))}.
\end{eqnarray*}

Here $b\rightarrow 0$ is a smoothing parameter (bandwidth),
$\Gamma(\cdot)$ is a standard gamma function and
  \vspace{-4mm}
\begin{eqnarray*}
\rho_b(x)&=& \left\{
\begin{array}{ll}
\rho_1(x) = \frac{x}{b}, &   \mbox{if}\qquad x\gs 2b,
\\
\rho_2(x) =\left(\frac{x}{2b}\right)^2+1, & \mbox{if}\qquad x\in
[0,2b).
\end{array}
\right.
\end{eqnarray*}

The support of the gamma kernel matches the support of the
probability density function to be estimated. For convenience let us
introduce two kernel functions
   \vspace{-4mm}
\begin{eqnarray*}
 K_{\rho_1(x),b}(t)&=\frac{t^{x/b-1}\exp(-t/b)}{b^{x/b}\Gamma(x/b)},
  \quad    &\mbox{if}\qquad x\gs 2b,\\
 K_{\rho_2(x),b}(t)&=\frac{t^{(\frac{x}{2b})^2}\exp(-t/b)}{b^{(\frac{x}{2b})^2+1}\Gamma((\frac{x}{2b})^2+1)},
 \quad &\mbox{if}\qquad x\in [0,2b).
\end{eqnarray*}

The estimate  $\hat{f}'(x)$ for density derivative $f'(x) =
df(x)/dx$ is usually taken as derivative of $\hat{f}(x).$ Hence, we
can write it as follows
\begin{gather}\label{2}
\hat{f}'(x) = \frac{1}{n}\sum_{i=1}^{n}K'_{\rho_b(x),b}(X_i),
\end{gather}
where
 \vspace{-4mm}
\begin{eqnarray*}
&&K'_{\rho_b(x),b}(t)= \\
&=&\left\{
\begin{array}{ll}
K'_{\rho_1(x),b}(t)=\frac{1}{b}K_{\rho_1(x),b}(t)L_1(t),
& \mbox{if}\quad x\gs 2b,\\
K'_{\rho_2(x),b}(t)=\frac{x}{2b^2}K_{\rho_2(x),b}(t)L_2(t),&
\mbox{if}\quad x\in [0,2b),
\end{array}
\right.
\end{eqnarray*}
with
 \vspace{-6mm}
\begin{eqnarray*}
L_1(t)=L_1(t;x)= \ln t - \ln b - \Psi(\rho_1(x)),\\
L_2(t)=L_2(t;x)= \ln t - \ln b - \Psi(\rho_2(x)).
\end{eqnarray*}
Here $\Psi(x)$  denotes  Digamma function (logderivative of gamma
function).

Now we get down to examine the properties of the derivative estimate
\eqref{2}. First of all, we investigate  the expectation $\mathsf
E(\hat{f}'(x)).$ It should be noted that each class of estimates has
nice properties only for a special class of densities. For instance,
the estimates proposed by Song Xi \cite{Chen:20} are matched to a
class of densities, satisfying the conditions: $f$ has a continuous
second derivative, and the integrals $\int_0^\infty f'^2(x)dx,$
$\int_0^\infty \{xf''(x)\}^2dx$ and $\int_0^\infty x^{-3/2}f(x)dx$
are finite. We intend to get analogous conditions for the density
derivative estimate  \eqref{2}.

{\bf Lemma 1.}(expectation) {\it If \ $b\rightarrow 0$ then the
leading term of the
 mathematical expectation expansion for the density derivative
estimate  \eqref{2} equals \vspace{-5mm}
\begin{eqnarray}
\mathsf E(\hat{f'}(x)) = \left\{
\begin{array}{ll}
\mathsf E K'_{\rho_1(x),b}(X_1), & \mbox{if}\quad x\gs 2b, \notag\\
\mathsf E K'_{\rho_2(x),b}(X_1), & \mbox{if}\quad x\in [0,2b),\notag
\end{array}
\right.
\end{eqnarray}
where }
\begin{gather*}
\begin{array}{ll}
&\mathsf E K'_{\rho_1(x),b}(X_1)=(1/b)\mathsf E
K_{\rho_1(x),b}(X_1)L_1(X_1;x)\\
&= f'(x) + b\left(\frac{1}{12x^2} f(x) +
\frac{1}{4}f''(x)\right)+o(b),\\
&\mathsf E K'_{\rho_2(x),b}(X_1)=(x/2b^2)\mathsf E K_{\rho_2(x),b}(X_1)L_2(X_1;x)\\
&=f'(x)\left(\frac{x}{2b}-\frac{b}{6x}\right)+f''(x)\left(\frac{7x}{48}+
\frac{x^2}{2b}\right)+o(b).
\end{array}
\end{gather*}
The proof of the Lemma 1 one can find in the Appendix. Note that
under fixed $b$  the estimate   $\hat{f'}(x)$ in the small area
$x\in [0,2b)$ near zero has a  bias, which  grows as $x\rightarrow
0$ . However, in the asymptotic case when $b\rightarrow 0$ the right
boundary of this area $x=2b$ decreases also to zero. Therefore, it
is interesting to know the bias limit when $x$ and $b$ converge to
zero simultaneously, i.e. when ratio $x/b$ tends to some constant
$\kappa$ when $b\rightarrow 0$. Then the second expectation of the
estimate  will differ very small from the true density derivative.
The leading term of bias expansion may be written as
  \vspace{-5mm}
\begin{eqnarray*}
Bias(\hat{f}'(x))&=&b\left(\frac{f(x)}{12x^2}+
\frac{f''(x)}{4}\right)+ o(b), \quad i\!f \;x/b\rightarrow \infty,\\
Bias(\hat{f}'(x))&=&f'(x)\left(\frac{3\kappa^2-6\kappa-1}{6\kappa}\right)\\
&+& b f''(x)\left(\frac{7\kappa}{48}+
\frac{\kappa^2}{2}\right)+o(b), \quad i\!f \;x/b\rightarrow \kappa.
\end{eqnarray*}
If $x=2b$ then $\kappa=2$ and the estimate  bias in the right
boundary of the small area near zero will differ from true density
derivative in $(1/12) f'(2b)$.

As a global performance of the density derivative estimate \eqref{2}
we select a mean integrated squared error ($MISE$), which, as is
known, equals to
  \vspace{-5mm}
\begin{eqnarray}\label{3}
MISE(\hat{f}'(x))&=&
\mathsf E\int\limits_0^\infty(f'(x)-\hat{f}'(x))^2dx \\
&=& \int\limits_0^\infty \big[Bias^2(\hat{f}'(x)) +
Var(\hat{f}'(x))\big]dx. \notag
\end{eqnarray}
As the right boundary $x=2b$ decreases with $n\rightarrow\infty,$
then the integral contribution to $MISE$ of the second part of the
bias in a small area near zero will be negligible. Hence, the
integral squared bias of the main area of support is important only.
Here it is
 \vspace{-4mm}
\begin{eqnarray*}
 IBias^2(\hat{f'}(x))
&=& \frac{b^2}{16} \int_{0}^\infty
\left(\frac{f(x)}{3x^2}+f''(x)\right)^2dx+o(b^2).
\end{eqnarray*}

Let us  proceed to calculate the variance of the derivative estimate.

{\bf Lemma 2.}(variance) {\it If  \ $b\rightarrow 0$ and $
nb^{3/2}\rightarrow \infty,$ then the leading term of
 variance expansion for density derivative estimate
\eqref{2} equals to}
 \vspace{-4mm}
 \begin{eqnarray*}&& Var(\hat{f}'(x))=
\\
&=& \frac{n^{-1}b^{-3/2}x^{-1/2}}{2\sqrt{\pi}}\left(\frac{f(x)}{2x}+b\left(\frac{f(x)}{4x^2}-\frac{f'(x)}{4x}\right)\right)+o(b)
\end{eqnarray*}
The proof of the Lemma 2 is in the Appendix.

The next problem is to calculate the mean squared error $MSE(x)$ in
accordance to the well known formula. Then
 \vspace{-4mm}
\begin{eqnarray*}
&&MSE\left(\hat{f}'(x)\right)=\frac{b^2}{16}\left(\frac{f(x)}{3x^2}\!\!\!\!+\!f''(x)\right)^2+\frac{n^{-1}b^{-3/2}x^{-1/2}}{2\sqrt{\pi}}\\
&\cdot&\left(\frac{f(x)}{2x}+b\left(\frac{f(x)}{4x^2}-\frac{f'(x)}{4x}\right)\right) + o(b^2),
\end{eqnarray*}
where
 \vspace{-5mm}
\begin{eqnarray*}
&&P(x)=\left(\frac{f(x)}{3x^2}+f''(x)\right)^2,
\end{eqnarray*}

If $P(x)\neq 0$, then minimization $MSE(x,b)= MSE\left(\hat{f}'(x)\right)$ in $b$ provides
an asymptotically optimal value of $b$
 \vspace{2mm}
\begin{gather}\label{4}
b_{opt}(x)=A(x)n^{-2/7},
\end{gather}
 where the so called initial coefficient $A(x)$ equals
 \vspace{-9mm}
\begin{eqnarray*}
&&A(x)=\left(\frac{3f(x)x^{-\frac{3}{2}}}{\sqrt{\pi}P(x)}\right)^{\frac{2}{7}}\!\!\!\!\!\!.
\end{eqnarray*}
The bandwidth  $b_{opt}(x)$ cannot be calculated directly because it
depends on the unknown true density $f(x).$ An algorithm for
evaluation $b_{opt}$ based on observations only   will be presented
in the next paper.

Substituting $b_{opt}(x)$ in $MSE\left(\hat{f}'(x)\right)$ leads to
the asymptotically optimal mean squared error of the estimate
$\hat{f}'_1(x)$ in each point $x$:
 \begin{eqnarray*}
&& MSE_{opt}(\hat{f}'(x))= \frac{A(x)^2P(x)^2}{16}n^{-\frac{4}{7}}+\frac{x^{-\frac{3}{2}}A(x)^{-\frac{3}{2}}}{4\sqrt{\pi}}n^{-\frac{4}{7}}.
\end{eqnarray*}

Now we proceed to global performance  \eqref{3}. We will receive the
integrated optimal bandwidth which doesn't depend on $x.$

{\bf Theorem} ($MISE$). {\it If $b\rightarrow 0$ and \ $
nb^{3/2}\rightarrow \infty,$ integrals
  \vspace{-7mm}
\begin{eqnarray*}
\int_{0}^\infty \!\!\!
\left(\frac{f(x)}{3x^2}+f''(x)\right)^2\!\!\!dx, \int_{0}^\infty
\!\!\!\!\! x^{-3/2} f(x)dx
\end{eqnarray*}
are finite and $\int_{0}^\infty P(x) dx\neq0$, then the leading term of a MISE expansion for the
density derivative estimate  $\hat{f}'(x)$ equals to
 \vspace{-5mm}
\begin{eqnarray} \label{5}
 &&MISE(\hat f'(x))=\frac{b^2}{16}\int_{0}^\infty
 \left(\frac{f(x)}{3x^2}+f''(x)\right)^2dx \nonumber
\\
&+&\!\!  \int_0^\infty \frac{n^{-1}b^{-3/2}x^{-3/2}}{4\sqrt{\pi}}\left(f(x)+\frac{b}{2}\left(\frac{f(x)}{x}-f'(x)\right)\right)dx
\\
&+&o(b^2 + n^{-1}b^{-3/2}). \nonumber
\end{eqnarray}
Minimization of \eqref{5} in $b$ leads to a global optimal bandwidth
 \vspace{-5mm}
\begin{eqnarray}\label{6}
 b_0 = \left(\frac{3\int_0^\infty x^{-3/2}f(x)dx}{\sqrt{\pi}\int_{0}^\infty
(\frac{f(x)}{3x^2}+f''(x))^2dx}\right)^{2/7}n^{-2/7},
\end{eqnarray}
whose substitution into \eqref{5} results to an optimal {$MISE$} }

The restrictions on the integrals in the Theorem are fulfilled, for
example, for the family of $\chi^2$-distributions with a number of
degrees of freedom $m\gs 3.$ For $m=3$ we receive Maxwell
distribution, which will be investigated as true distribution in
simulation below.

From expression for $MISE_{opt}$ it follows that nonparametric
estimate  \eqref{2} converges in mean square to true density
derivative with the rate $O(n^{-4/7}).$  This rate is certainly less
than the rate of convergence for the density $O(n^{-4/5})$, because the estimation of derivatives is more complex
than the estimation of the densities. A similar decrease in the rate of convergence
for the derivatives compared with the densities was observed in the
use of Gaussian kernel functions on the whole line.

\section{Simulation results}
In the simulation experiment we select the density of Maxwell
distribution with parameter $\sigma=1$ as the true density to be
estimated:
  \vspace{-8mm}
\begin{eqnarray*}
f_M(x)=\frac{\sqrt{2}x^2\exp(-\frac{x^2}{2\sigma^2})}{\sigma^3\sqrt{\pi}}.
\end{eqnarray*}
We need two derivatives of it for computation integrals in the
optimal bandwidth $b_0$ \eqref{6}
  \vspace{-6mm}
\begin{eqnarray*}
f'_M(x)=-\frac{\sqrt{2}x\exp(-\frac{x^2}{2\sigma^2})(x^2-2\sigma^2)}{\sigma^5\sqrt{\pi}},
\end{eqnarray*}
  \vspace{-6mm}
\begin{eqnarray*}
f''_M(x)=\frac{\sqrt{2}\exp(-\frac{x^2}{2\sigma^2})(2\sigma^4-5\sigma^2x^2+x^4)}{\sigma^7\sqrt{\pi}}.
\end{eqnarray*}

Sample sizes are $n=200$ and $n=2000$.  For comparison, the  values
of bandwidths were determined by three methods. The first calculates
bandwidth from the formula \eqref{6}. It has two integrals where we
have to substitute $f_M(x)$ and its derivatives instead of $f(x)$
with corresponding derivatives. For $\sigma=1$ and $n=2000$ the
values of the integrals are:
\\in numerator
  \vspace{-8mm}
\begin{eqnarray*}\left(\frac{3}{\sqrt{\pi}}\int_0^\infty x^{-3/2}f_M(x)dx\right)^{2/7}=
1.099, \qquad n^{-2/7}=0.114;
\end{eqnarray*}
\\in denominator
  \vspace{-4mm}
\begin{eqnarray*}\left(\int_{0}^\infty
\left(\frac{f_M(x)}{3x^2}+f''_M(x)\right)^2dx\right)^{2/7}=1.247.
\end{eqnarray*}
Combining these data together, we obtain $b_1 = 0.1004$.
\\The second
bandwidth is a solution of equation \eqref{11}, where the integral
coefficients were calculated numerically. The coefficient of $b$ is
 \vspace{-8mm}
\begin{eqnarray*}
\frac{1}{8}\left(\int_{0}^\infty
(\frac{f_M(x)}{3x^2}+f''_M(x))^2dx\right)^{2}=0.270.
\end{eqnarray*}
The coefficients of  $b^{-5/2}$ and $b^{-3/2}$  are, respectively,
 \vspace{-6mm}
\begin{eqnarray*}
&&\frac{3n^{-1}}{8\sqrt{\pi}}\!\int_0^\infty\!\!\!x^{-\frac{3}{2}}f_M(x)dx=8.69\cdot10^{5},\\\nonumber
&&\frac{n^{-1}}{16\sqrt{\pi}}\!\int_0^\infty\!\!\!x^{-\frac{3}{2}}\left(\frac{f_M(x)}{x}-f_M'(x)\right)dx=-2.69\cdot10^{5}.
\end{eqnarray*}
Substitution all of  them in \eqref{11} yields the transcendental
equation
 \vspace{-8mm}
\begin{eqnarray*}
0.270b+8.69\cdot10^{5}b^{-5/2}-2.69\cdot10^{5}b^{-3/2}=0,
\end{eqnarray*}
which can be solved by numerical methods. Solution of it provides
 $b_2=0.1013$.

The third $b_3$ was taken from the paper of Song Xi Chen (2000),
where he found an optimal in mean square sense bandwidth
 \vspace{-6mm}
\begin{eqnarray*}
b_3 = \left(\frac{V}{\beta}\right)^{2/5}n^{-2/5},
\end{eqnarray*}
 where
  \vspace{-6mm}
\begin{eqnarray*}
V = \frac{1}{2\sqrt{\pi}}\int\limits_0^{\infty} x^{-1/2}f_M(x)dx,
\quad \beta = \ilt_0^\infty (xf_M''(x))^2dx\neq0.
\end{eqnarray*}
In our case it is equals to $b_3=0.0175$. In this case one might
think that  if the estimate of the derivative of density is
constructed as a derivative of the density estimate, then the best
bandwidth for
 the density will be good for its derivative. However,
this is not the case, as evidenced by the experimental results.
\begin{figure}
\begin{center}
\includegraphics[width=8.4cm]{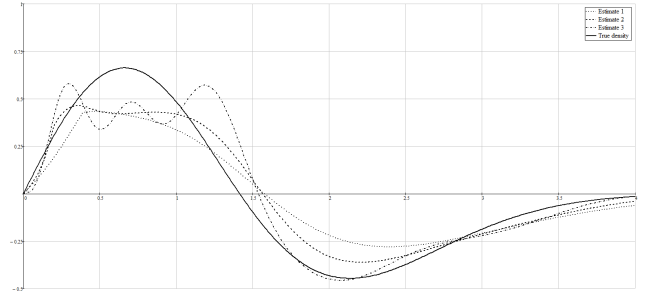}
\caption{Nonparametric estimates of Maxwell density derivative
function for n=200. The $f'_M(x)$ (solid line), estimate 1 b1=0.194
(dotted line), estimate 2 b2=0.197 (dashed  line), estimate 3
b3=0.0175 (dash-dotted line).} \label{fig:2}
\end{center}
\end{figure}
\begin{figure}
\begin{center}
\includegraphics[width=8.4cm]{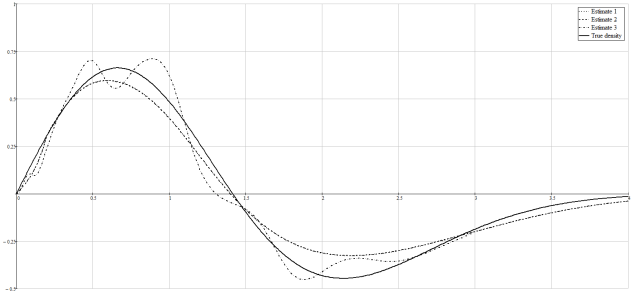}
\caption{Nonparametric estimates of Maxwell density derivative
function for n=2000. The $f'_M(x)$ (solid line), estimate 1 b1=0.1004
(dotted line), estimate 2 b2=0.1013 (dashed  line), estimate 3
b3=0.0175 (dash-dotted line).} \label{fig:1}
\end{center}
\end{figure}

From Fig. 1 and Fig. 2 it is visible that the estimate 1 and 2 areclose to the desire density derivative (solid
line). The estimate 2 is quite
smooth compared with the estimate 3 using the optimal bandwidth
$b_3$ for the density. This is confirmed by the numerical evaluation
of the squared integral deviation from the true derivative curve, as
it is shown in the table.

\begin{table}[htbp]
\begin{center}
\caption{} {\normalsize  Deviation of estimates }
\end{center}
\vspace{2mm}
\begin{center}
\begin{tabular}{|c|c|c|c|} \hline
Bandwidth &  $b_1$  &  $b_2$  & $b_3$ \\
\hline Value& 0.203 & 0.146   &  0.017\\
\hline Deviation& 0.0426 & 0.0382 & 0.0450 \\
\hline
\end{tabular}
\end{center}
\end{table}

The estimate 2 with $b_2$, when we solve a transcendental equation,
provides the best result. If there are multiple roots of the
transcendental equation, we choose the root with the lowest value of
$MISE$.

\section{Conclusions}
We have developed a method of nonparametric density derivative
estimation on the positive semi-axis, which is supposed to be applied
to nonlinear problems of signal selection with unknown
characteristics from the mixture with noise. Such problems are
arisen in the theory of nonparametric estimation of signals with
unknown distribution, where there is an equation of optimal
filtering, containing statistics in the form of the logarithmic
derivative of the density. In multiplicative observation models with
positive signals the logarithmic derivative has to be reconstructed
from observations. Since the logarithmic derivative contains a
derivative of the unknown density, the presented method for
estimating the derivative is relevant.

This method is expected to be extended to dependent variables.
In addition, since the optimal bandwidth depends on the unknown
density, it is necessary to build its data-based estimate and thus
to create an automatic technique of nonparametric signal
estimation.

\bibliography{reference}

\begin{thebibliography}{9}
\providecommand{\natexlab}[1]{#1}
\providecommand{\url}[1]{\texttt{#1}}
\providecommand{\urlprefix}{URL }
\expandafter\ifx\csname urlstyle\endcsname\relax
  \providecommand{\doi}[1]{doi:\discretionary{}{}{}#1}\else
  \providecommand{\doi}{doi:\discretionary{}{}{}\begingroup
  \urlstyle{rm}\Url}\fi

\bibitem[{Chen(2000)}]{Chen:20}
Song Xi Chen. (2000).
\newblock Probability density function estimation using gamma kernels.
\newblock \emph{Ann. Inst. Statist. Math., Australia.}

\bibitem[{Jones(1995)}]{Jones:95}
Wand, M. and Jones, M.C. (1995).
\newblock \emph{Kernel Smoothing.}
\newblock \emph{Chapman and Hall, London.}

\bibitem[{Brown and Chen(1999)}]{Brown:99}
Brown, B, M. and Chen, S. X. (1999).
\newblock Beta-Bernstein smoothing for regression curves with compact support.
\newblock \emph{Scand. J. Statist.}

\bibitem[{ Dobrovidov et~al. Dobrovidov, Koshkin, and Vasiliev}]{Dobrovidov:12}
Dobrovidov, A., Koshkin, G., and Vasiliev, V. (2012).
\newblock Non-Parametric state space models.
\newblock \emph{Kendrick press}

\end{thebibliography}
\vspace{5mm}

Appendix:  PROOFS OF THE STATEMENTS.

{\bf Proof of Lemma 1.}

We start by writing the expectation of the estimate \eqref{2}
on the both intervals of the support.\\
1) For the interval $x\gs2b$
 \vspace{-5mm}
\begin{eqnarray}\label{7}
&& \mathsf{E}(\hat{f}'(x))=
\int_{0}^{\infty}\frac{1}{b}K_{\rho_1,b}(y)L_1(y;x)f(y)dy \nonumber
\\
&=& \mathsf{E}\left(b^{-1}(\ln \xi_x - \ln b - \Psi(\rho_1))\cdot
f(\xi_x)\right) \nonumber
\\
&=& \frac{1}{b}[\mathsf{E}(\ln \xi_x \cdot f(\xi_x)) - (\ln b +
\Psi(\rho_1))\mathsf{E}( f(\xi_x))],
\end{eqnarray}
where $\xi_x$ is a $\text {Gamma}(\rho_1,b)$ random variable. From
the standard theory of gamma distribution it is known that for this
random variable mean is $\mu_{x} = \mathsf{E}(\xi_x)= \rho_1 b$ and
variance is $Var_{1}(\xi_x)=\rho_1 b^2$.
\\
2) Regarding the interval $x\in[0,2b)$ we get in a similar way:
 \vspace{-5mm}
\begin{eqnarray}\label{8}
&& \mathsf{E}(\hat{f'}(x))= \int_{0}^{\infty}K_{\rho_2,b}(y)L_2(y;x)\frac{x}{2b^2}f(y)dy \nonumber
\\
&=& \frac{x}{2b^2} \mathsf{E}((\ln \eta_x - \ln b -
\Psi(\rho_2))\cdot f(\eta_x)) \nonumber
\\
&=& \frac{x}{2b^2} [\mathsf{E}(f(\eta_x)\ln \eta_x) - (\ln b +
\Psi(\rho_2)\cdot \mathsf{E}(f(\eta_x))],
\end{eqnarray}
 where $\eta_x$ is a $\text {Gamma}(\rho_2,b)$ random variable with mean
$\mu_{x} = \mathsf{E}(\eta_x)= \rho_{2} b$ and variation
$Var_{2}(\eta_x)=\rho_{2} b^2$. Up to a factor, it is the same as in
\eqref{7}. Then  we will make the Taylor series expansion at a point
$\mu_x$ for general $\rho$  and then substitute in the appropriate
cases $\rho_1$ or $\rho_2$. Consider the first term in
\eqref{7},\eqref{8}:
 \vspace{-5mm}
\begin{eqnarray*}
&&\mathsf{E}(f(\xi_x)\ln \xi_x) = \mathsf{E}(f(\mu_x)\ln \mu_x)
\\
&+& \mathsf{E}((f(\mu_x)\ln\mu_x)'(\xi_x - \mu_x))
\\
&+& \frac{1}{2}\mathsf{E}((f(\mu_x)\ln\mu_x)''(\xi_x -
\mu_x)^2)+o(b)\\
&=&f(\mu_x)\ln\mu_x+\frac{Var(\xi_x)}{2}
\\
&\cdot&\left(\frac{2f'(\mu_x)}{\mu_x}+f''(\mu_x)\ln(\mu_x)-\frac{f(\mu_x)}{\mu_x^2}\right)+
o(b).
\end{eqnarray*}

We substitute the mean and variance by their values from a gamma
 distribution:
  \vspace{-5mm}
\begin{eqnarray*}
&& f(\rho b)(\ln(\rho b)\!-\frac{1}{2\rho})\!+\!f'(\rho b)b\!+\!f''(\rho b)\ln(\rho b)\frac{(\rho b)^2}{2}+o(b)\\
&=&f(x)\!\left(\ln(\rho b)\!-\frac{1}{2\rho}\right)\!+\!f'(x)\left(\!\!\left(\ln(\rho b)-\frac{b}{2}\right)(x-\rho b)+b\!\right)\\
&+&f''(x)\lt(\left(\ln(\rho b)-\frac{1}{2\rho}\right)\frac{(x-\rho b)^2}{2}+b(x-\rho b) \rt.\\
&+& \lt. \ln(\rho b)\frac{(\rho b)^2}{2}\rt)+o(b).
\end{eqnarray*}

The second term in \eqref{7}\eqref{8} can be represented just like
above
 \vspace{-5mm}
\begin{eqnarray*}&&\mathsf{E}(f(\xi_x))\!=\!f(x)+f'(x)(x-\rho b)\!+\!f''(x)\!\left(\!\frac{\rho b^2}{2}\!+\!\frac{(x-\rho)^2}{2}\!\right)
\\
&+&o(b).
\end{eqnarray*}

Then, combining all items in the square brackets in
\eqref{7},\eqref{8}, we can get
 \vspace{-5mm}
\begin{eqnarray*}&&[f(x)\left(\ln(\rho b)-\frac{1}{2\rho}-(\ln b +\Psi(\rho))\right)\\
&+&\!f'(x)\!\left(\!\!\left(\ln(\rho b)-\frac{1}{2\rho}\right)\!(x-\rho b)+\!b-\!(\ln b +\Psi(\rho))(x\!-\rho b)\!\right)\\
&+&\!f''(x)(\left(\ln(\rho b)-\frac{1}{2\rho}\right)\!\frac{(x-\rho b)^2}{2}\!+\!b(x-\rho b)+\ln(\rho b)\frac{(\rho b)^2}{2}\\
&-&(\ln b +\Psi(\rho))(\frac{\rho b^2}{2}+\frac{(x-\rho
b)^2}{2}))+o(b)].
\end{eqnarray*}

Using the approximation of the Digamma function when $\rho\rightarrow\infty$
 \vspace{-4mm}
\begin{eqnarray*}\Psi(\rho) &=& \ln \rho - \frac{1}{2\rho}-\frac{1}{12\rho^2}+\frac{1}{120\rho^4}-
\frac{1}{252\rho^6}+O\left(\frac{1}{\rho^8}\right),
\end{eqnarray*}

we receive the expression in square brackets of \eqref{8} in  $\rho$
 \vspace{-5mm}
\begin{eqnarray*}
&&[f(x)\left(\frac{1}{12\rho^2}\right)+f'(x)\left(b+\frac{x-\rho b}{12\rho^2}\right)\\
&+&\!f''(x)\!\left(\!\left(\!\frac{1}{2\rho}\!+\!\frac{1}{12\rho^2}\right)\frac{\rho
b^2}{2}\!+\!\frac{(x-\rho b)^2}{24\rho^2}\!+b(x-\rho
b)\!\right)\!+\!o(b)].
\end{eqnarray*}
Now for cases 1) and 2) let us  substitute $\rho_{1}$ and $\rho_{2}$
instead of $\rho$.
\\
For 1):
 \vspace{-5mm}
\begin{eqnarray*}\mathsf E(\hat{f}'(x))
&=& \frac{1}{b}\left(f(x)\left(\frac{b^2}{12x^2}\right) + b f'(x) +
\frac{b^2}{4}f''(x)+ o(b)\right)
\\
&=&f'(x) + b\left(\frac{1}{12x^2} f(x) +
\frac{1}{4}f''(x)\right)+o(b).
\end{eqnarray*}
For 2) we will use the fact that as $b\rightarrow 0$ then
 \vspace{-5mm}
\begin{eqnarray*}&&\frac{1}{\rho}=\frac{4b^2}{x^2(1+\frac{4b^2}{x^2})}=
\frac{4b^2}{x^2}+o(b^2),
\end{eqnarray*}
  and
 \vspace{-8mm}
\begin{eqnarray*}
&&\mathsf
E(\hat{f'}(x))=f'(x)\left(\frac{x}{2b}-\frac{b}{6x}\right)+f''(x)\left(\frac{7x}{48}+
\frac{x^2}{2b}\right)+o(b).\\ &&\Box
\end{eqnarray*}

{\bf Proof of Lemma 2.}

We start with variance for $x\gs 2b$.
 \vspace{-8mm}
\begin{eqnarray}\label{9}
&& Var(\hat{f}'(x))=\frac{1}{n}Var(K'_{\rho_1,b}(x)) \nonumber
\\
&=&\frac{1}{n}\big(\mathsf E(K'^2_{\rho_1,b}(x))-\mathsf
E^2(K'_{\rho_1,b}(x))\big).
\end{eqnarray}
The second term of the right-hand side of \eqref{9} is the same as
in \eqref{7}.  So we can  write immediately
 \vspace{-5mm}
\begin{eqnarray*}&&\mathsf E^2(K'_{\rho_1,b}(x))=\left(f(x)\frac{b}{12x^2}+
f'(x)+f''(x)\frac{b}{4}+o(b)\right)^2.
\end{eqnarray*}
The first term of the right-hand side of (\ref{9}) can be
represented by
 \vspace{-5mm}
\begin{eqnarray*}&&\mathsf E(K'^2_{\rho_1,b}(x))=\int_0^\infty K'^2_{\rho_1,b}(y) f(y)dy\\
&=&\int_0^\infty\frac{y^{\frac{2x}{b}-2}\exp(\frac{-2y}{b})}{b^{\frac{2x}{b}}\Gamma^2(\frac{x}{b})}
\left(\frac{L_1(y;x)}{b}\right)^2 f(y)dy.
\end{eqnarray*}
Using the property of the gamma function
$\Gamma^2(\frac{x}{b}+1)=(\frac{x}{b})^2\Gamma^2(\frac{x}{b}),$ we
get
 \vspace{-10mm}
\begin{eqnarray*}
\\
&&
\int_0^{\infty}\frac{y^{\frac{2x}{b}-2}\exp(\frac{-2y}{b})}{b^{\frac{2x}{b}}(\frac{b}{x})^2\Gamma^2(\frac{x}{b}+1)}\cdot
\left(\frac{L_1(y;x)}{b}\right)^2 f(y)dy
\\
&=&\int_0^{\infty}\frac{b^{-5}x^2\Gamma(\frac{2x}{b}-1)}{2^{\frac{2x}{b}-2}\Gamma^2(\frac{x}{b}+1)}\frac{(2y)^{\frac{2x}{b}-2}\exp(\frac{-2y}{b})}{b^{\frac{2x}{b}-1}\Gamma^(\frac{2x}{b}-1)}
\\
&\cdot & L_1(y;x)^2 f(y)dy.
\end{eqnarray*}
Denoting $B_b(x) =
\frac{b^{-5}x^2\Gamma(\frac{2x}{b}-1)}{2^{\frac{2x}{b}-2}\Gamma^2(\frac{x}{b}+1)},$
it can be written shorter
 \vspace{-5mm}
\begin{eqnarray}\label{10}
&=&\int_0^\infty B_b(x)K_{\frac{2x}{b}-1,b}(y)L_1(t)^2 f(y)dy \nonumber \\
&=& B_b(x)\mathsf E(L_1(\eta_x,x)^2f(\eta_x)),
\end{eqnarray}
where $\eta_x$ is a Gamma$(\frac{2x}{b}-1,b)$  random variable with
a mean  $\mu_x = \mathsf E(\eta_x)= 2x-b$ and a  variance
$Var(\eta_x)=2xb-b^2$. Let $R(z) =
\sqrt{2\pi}\exp(-z)z^{z+1/2}/\Gamma(z+1)$ for $z\gs0$. So we can
express gamma function as
 \vspace{-6mm}
\begin{eqnarray*}\Gamma^2(\frac{x}{b}+1)&=&\left(\frac{\sqrt{2\pi}\exp(\frac{-x}{b})
(\frac{x}{b})^{\frac{x}{b}+1/2}}{R(\frac{x}{b})}\right)^2.
\end{eqnarray*}
Using the properties of the gamma function
 \vspace{-4mm}
\begin{eqnarray*}\Gamma(\frac{2x}{b}-1)=\frac{\Gamma(\frac{2x}{b}+1)}{\frac{2x}{b}(\frac{2x}{b}-1)}=
\frac{\sqrt{2\pi}\exp(-\frac{2x}{b})(\frac{2x}{b})^{\frac{2x}{b}+\frac{1}{2}}}{\frac{2x}{b}(\frac{2x}{b}-1)R(\frac{2x}{b})},
\end{eqnarray*}
we obtain
 \vspace{-5mm}
\begin{eqnarray*}&&B_b(x)=\frac{b^{-5}x^2}{2^{\frac{2x}{b}-2}}\cdot \frac{1}{\frac{2x}{b}(\frac{2x}{b}-1)}\frac{\sqrt{2\pi}\exp(-\frac{2x}{b})
(\frac{2x}{b})^{\frac{2x}{b}+\frac{1}{2}}}{R(\frac{2x}{b})}\\
&\cdot &\frac{R^2(\frac{x}{b})}{(\sqrt{2\pi})^2\exp(-\frac{2x}{b})
(\frac{x}{b})^{\frac{2x}{b}+1}}
=
\frac{b^{-\frac{5}{2}}x^{-\frac{1}{2}}R^2(\frac{x}{b})}{\sqrt{\pi}R(\frac{2x}{b})(1-\frac{b}{2x})}.
\end{eqnarray*}
According to Lemma 3 of \cite{Brown:99}, $R(z)$ is
increasing function which converges to 1 as $z\rightarrow\infty$ and
$R(z)<1$ for any $z>0$.

 Then
 \vspace{-10mm}
\begin{eqnarray*}
B_b(x)&=& \left\{ \begin{array}{ll}
\frac{b^{-5/2}x^{-1/2}}{2\sqrt{\pi}}, &   \mbox{if}\qquad
\frac{x}{b}\rightarrow\infty,
\\
\frac{b^{-3}k^2\Gamma(2k-1)}{2^{2k-2}\Gamma^2(k+1)}, &
\mbox{if}\qquad \frac{x}{b}\rightarrow k,
\end{array}
\right.
\end{eqnarray*}
Let us denote $G = G(x,b)=\ln b + \Psi(\frac{x}{b}) = \ln x  -
\frac{b}{2x}-\frac{b^2}{12x^2}+o(b^2)$. Now we must find the second
factor in \eqref{10}
 \vspace{-7mm}
\begin{eqnarray*}
&&\mathsf E(L_1(\eta_x;x)^2f(\eta_x))=\mathsf E\left(\big(\ln \eta_x - \ln b - \Psi(\frac{x}{b})\big)^2 f(\eta_x)\right)\\
&=&\mathsf E(f(\eta_x)\ln^2 \eta_x) -  2G\mathsf{E}(f(\eta_x)\ln
\eta_x)+G^2\mathsf{E}(f(\eta_x)).
\end{eqnarray*}
If $b/x\rightarrow 0,$ we can get Taylor series
 \vspace{-8mm}
\begin{eqnarray*}
\ln(2x-b)&=&\ln(2x(1-\frac{b}{2x}))=\ln(2x)+ \ln(1-\frac{b}{2x})
\\
&=& \ln(2x)- \frac{b}{2x}+ o(b),
\end{eqnarray*}
 \vspace{-12mm}
\begin{eqnarray*}\frac{1}{2x-b}=\frac{1}{2x(1-\frac{b}{2x})}=
\frac{1}{2x}+\frac{b}{4x^2}+o(b).
\end{eqnarray*}
Substituting them in the expression above, we obtain
 \vspace{-7mm}
\begin{eqnarray*}&&
\mathsf E(f(\eta_x)\ln^2 \eta_x)=f(\mu_x)\ln^2 \mu_x + (f(\eta_x)\ln^2 \eta_x)''\frac{Var(\eta_x)}{2}\\
&=& f(x)\ln^2(x)+b\Big(f'(x)\left(\ln(x)-\frac{\ln^2(x)}{2}\right)\\
&+&\frac{f(x)}{2x}\left(1-3\ln(x)\right)+f''(x)\frac{x\ln^2(x)}{4}\Big)\\
&+&b^2\Big(f(x)\left(\frac{3}{4x^2}-\frac{\ln(x)}{2x^2}\right)+f'(x)\frac{3}{4x}\left(\ln(x)-1\right)\\
&+&f''(x)\left(\frac{\ln^2(x)}{8}-\frac{3\ln(x)}{4}\right)-f'''(x)\frac{x\ln^2(x)}{8}\Big)+o(b).
\end{eqnarray*}
Similarly,
  \vspace{-7mm}
\begin{eqnarray*}&&
\mathsf{E}(f(\eta_x)\ln \eta_x)=f(\mu_x)\ln \mu_x + (f(\eta_x)\ln \eta_x)''\frac{Var(\eta_x)}{2}\\
&=&f(x)\ln(x)+b\Big(-f(x)\frac{3}{4x}+f'(x)\frac{\left(1-\ln(x)\right)}{2}\\
&+&f''(x)\frac{x\ln(x)}{4}\Big)+b^2\Big(-\frac{f(x)}{4x^2}f'(x)\frac{3}{8x}\\
&+&f''(x)\left(\frac{\ln(x)}{8}-\frac{3}{8}\right)-f'''(x)\frac{x\ln(x)}{8}\Big)+o(b),
\end{eqnarray*}
  \vspace{-11mm}
\begin{eqnarray*}&&
\mathsf{E}(f(\eta_x))=f(\mu_x)+f''(\mu_x)\frac{Var(\xi_x)}{2}
\\
&=& f(x)-b\left(\frac{xf''(x)}{4}-\frac{f'(x)}{2}\right)+b^2\left(\frac{f''(x)}{8}-\frac{xf'''(x)}{8}\right)+o(b).
\end{eqnarray*}

Collecting all the terms, we obtain an expression
 \vspace{-7mm}
\begin{eqnarray*}&&
 \mathsf{E}((\ln \eta_x - \ln b - \Psi(x/b))^2 f(\eta_x))
\\
&=& b\frac{f(x)}{2x}+b^2\left(\frac{f(x)}{4x^2}-\frac{f'(x)}{4x}\right)+o(b).
\end{eqnarray*}
Hence, as $b/x\rightarrow 0$, the variance is
 \vspace{-7mm}
\begin{eqnarray*}&&
Var(\hat{f}'_2(x))=
\\
&=& \frac{n^{-1}b^{-3/2}x^{-1/2}}{2\sqrt{\pi}}\left(\frac{f(x)}{2x}+b\left(\frac{f(x)}{4x^2}-\frac{f'(x)}{4x}\right)\right).\quad\Box
\end{eqnarray*}
\begin{multline*}
\end{multline*}
{\bf Proof of Theorem.}

Once the variance of the estimate  \eqref{2} is calculated, we
simply use the expression \eqref{3} to obtain formula \eqref{5}.
Differentiation of the last expression in $b$ leads to equation
 \vspace{-6mm}
\begin{eqnarray}\label{11}
&&\!\frac{b}{8}\int_{0}^\infty\!\!\!
\left(\frac{f(x)}{3x^2}+f''(x)\right)^2\!\!dx\!-\!
\frac{3n^{-1}b^{-\frac{5}{2}}}{8\sqrt{\pi}}\!\int_0^\infty\!\!\!x^{-\frac{3}{2}}f(x)dx\\\nonumber
&+&\frac{n^{-1}b^{-\frac{3}{2}}}{16\sqrt{\pi}}\!\int_0^\infty\!\!\!x^{-\frac{3}{2}}\left(\frac{f(x)}{x}-f'(x)\right)dx=0
\end{eqnarray}
If we  neglect the term with $b^{-3/2}$ as compared to the term with
$b^{-5/2}$,  the equation becomes simpler
and  its solution is  equal to the optimal global bandwidth $b_0$.
\; $\Box$
\\Nevertheless, the use of equation \eqref{11} is also useful, because
its numerical solution gives $b_0'$  which, as shown in simulation,
yields a little better quality with respect to the case with $b_0$.


\end{document}